\documentclass{amsart}

\usepackage{amsmath,amsthm}
\usepackage{amsfonts,amssymb}

\newtheorem{thm}{Theorem}[section]
\newtheorem{cor}[thm]{Corollary}
\newtheorem{lem}[thm]{Lemma}
\newtheorem{prop}[thm]{Proposition}

\theoremstyle{remark}

 \def\ab{{\mathbf a}}
 \def\xb{{\mathbf x}}
 \def\yb{{\mathbf y}}
 \def\Kb{{\mathbf K}}
 \def\Sb{{\mathbf S}}
 \def\CL{{\mathcal L}}
 \def\CO{{\mathcal O}}
 \def\CS{{\mathcal S}}
 \def\CV{{\mathcal V}}
 
 \def\PP{{\mathbb P}}
 \def\QQ{{\mathbb Q}}
 \def\RR{{\mathbb R}}
 
     \def\Pad{\operatorname{Pad}}
     \def\codim{\operatorname{codim}}
     \def\sspan{\operatorname{span}}
     \def\LT{\operatorname{LT}}
\newcommand{\wt}{\widetilde}

\begin{document}

\title[Bivariate Lagrange Interpolation at Padua points]
{Bivariate Lagrange interpolation at the Padua points:  
the ideal theory approach\thanks{Work 
supported by the ``ex-$60\%$'' funds of the University of Padua and  
by the INdAM GNCS (Italian National Group for Scientific Computing).}}
\author{Len Bos} 
\address{Department of Mathematics and Statistics\\ University of Calgary\\
Calgary, Alberta, Canada T2N1N4}
\email{lpbos@math.ucalgary.ca}
\author{Stefano De Marchi}
\address{Department of Computer Science\\ University of Verona\\ 37134 Verona, Italy}
\email{demarchi@sci.univr.it}
\author{Marco Vianello}
\address{Department of Pure and Applied Mathematics\\ University of Padua\\35131 Padova,  Italy} 
\email{ marcov@math.unipd.it}
\author{Yuan Xu}
\address{  Department of Mathematics\\ University of Oregon\\
    Eugene, Oregon 97403-1222 \\ USA.}
\email{yuan@math.uoregon.edu}
 
\date{\today}
\keywords{bivariate Lagrange interpolation, Padua points, polynomial ideal}
\subjclass{41A05, 41A10}

\begin{abstract}
Padua points is a family of points on the square $[-1,1]^2$ given by explicit
formulas that admits unique Lagrange interpolation by bivariate polynomials.
The interpolation polynomials and cubature formulas based on the Padua
points are studied from an ideal theoretic point of view, which leads to the 
discovery of a compact formula for the interpolation polynomials. The $L^p$ 
convergence of the interpolation polynomials is also studied. 
\end{abstract}

\maketitle

%\begin{AMS}
%{\rm ???.} 
%\end{AMS}

\section{Introduction.}
\setcounter{equation}{0}

For polynomial interpolation in one variable, the Chebyshev points (zeros
of Chebyshev polynomial) are optimal in many ways.  Let $L_n f$ be the
$n$-th Lagrange interpolation polynomial based on the Chebyshev points.
The Lebesgue constant $\|L_n\|$ in the uniform norm on $[-1,1]$ grows in 
the order of $\CO(\log n)$, which is the minimal rate of growth of any 
projection 
operator from $C[-1,1]$ onto the space of polynomials of degree at most $n$. 
The Chebyshev points are also the knots of the Gaussian quadrature formula 
for the  weight function $1/\sqrt{1-x^2}$ on $[-1,1]$. 

In the case of two variables, it is often difficult to give explicit point 
sets that are unisolvent for polynomial interpolation. The Padua points is a 
set of points that may be considered as an analogue of the Chebyshev points
on the square $[-1,1]^2$.  Let $\Pi_n^2$ denote the space of  polynomials 
of degree at most $n$ in two variables. It is well known that $\dim \Pi_n^2 
= (n+2)(n+1)/2$.  Let $V$ be a set of points in $\RR^2$ and assume that 
the cardinality of $V$ is $|V| = \dim \Pi_n^2$. The set  $V$ is said to be 
unisolvent if for any given function $f$ defined on $\RR^2$, there is a 
unique polynomial $P$ such that $P(x) = f(x)$ on $V$. For each $n \ge 0$,
the Padua points are defined by  
\begin{equation} \label{padua1} 
  \Pad_n=\{ \xb_{k,j} =  (\xi_k,\eta_j), \quad 0 \leq k\leq n, \quad 
         1\leq j\leq \lfloor \tfrac{n}{2} \rfloor +1 \}, 
\end{equation}
where 
\begin{equation} \label{padua2}
\xi_k=\cos{\frac{k \pi}{n}}, \qquad \eta_j=\left\{ \begin{array} {lr}
   \cos{\frac{2j-1}{n+1}\pi}, & \;k \;\mbox{even}\\  &  \\
\cos{\frac{2j-2}{n+1}\pi},  & \;k \;\mbox{odd} 
\end{array}
\right.
\end{equation}
It is easy to verify that the cardinality of $\Pad_n$ is equal to the dimension
of $\Pi_n^2$. These points were introduced heuristically in  \cite{CDMV05} 
(only for even degrees) and proved to be unisolvent in \cite{BCDMVX}. The 
Lebesgue constant of the Lagrange interpolation based on the Padua points 
in the uniform norm 
grows like $\CO((\log n)^2)$ and there is a cubature formula of degree 
$2n-1$ based on these points for the weight function
$1/ ( \sqrt{1-x^2}\sqrt{1-y^2} )$ on $[-1,1]^2$ (\cite{BCDMVX}).  A family of 
points studied in \cite{X96} also possesses  similar properties 
(\cite{BCDMV05,BDM,X96}), but the cardinality of the set is not equal to 
$\dim \Pi_n^2$ and the interpolation polynomial belongs to a subspace of 
$\Pi_n^2$. 

The study in \cite{BCDMVX} starts from the fact that the Padua 
points there 
lie on a single generating curve, 
$$
    \gamma_n(t) = (\cos (n t), \cos ((n+1)t), \qquad 0 \le t \le \pi.
$$
In fact, the points are exactly the distinct points among 
$\{\gamma_n(\frac{k\pi}{n(n+1)})$, $0 \le k \le n(n+1)\}$. 
In the present paper, we study the interpolation on the Padua points using
an ideal theoretic approach, which considers the points as the variety of a
polynomial ideal and it treats the cubature formula based on the 
interpolation
points simultaneously. The advantage of this approach is that it casts the 
result on Padua points into a general theoretic framework.  The study in 
\cite{X96} used such an approach. In particular, it shows how the compact
formula of the interpolation polynomial arises naturally. 

The order of the Lebesgue constant of the interpolation is found in 
\cite{BCDMVX}, which solves the problem about the convergence of the 
interpolation polynomials on the Padua points in the uniform norm. In the 
present paper, we show the convergence in a weighted $L^p$ norm. 

The paper is organized as follows. In the following section, we review the
necessary background and show that the Padua points are unisolvent. 
The explicit formula of the Lagrange interpolation polynomials is derived
in Section 3, and the $L^p$ convergence of the interpolation polynomials
is proved in Section 4. 
 
\section{Polynomial ideal and Padua points} 
\setcounter{equation}{0}

\subsection{Polynomial ideals, interpolation, and cubature formulas}

We first recall some notation and results about polynomial ideals and their
relation to polynomial interpolation and cubature formulas. To keep the 
notation simple we shall restrict to the case of two variables, even though
all results in this subsection hold for more than two variables. 

Let $I$ be a polynomial ideal in $\RR[\xb]$ with $\xb =(x_1,x_2)$, the 
ring of all polynomials in two variables. If there are polynomials $f_1,
\ldots, f_r$ in $\RR[\xb]$ such that every $f \in I$ can be written as 
$f = a_1 f_1 + \ldots + a_r f_r$, $a_j \in \RR[\xb]$, then we say that $I$ is 
generated by the basis $\{f_1,\ldots,f_r\}$ and we write $I = \langle f_1,
\ldots, f_r \rangle$.  Fix a monomial order, say the lexicographical order,  
and let $\LT(f)$ denote the leading term of the polynomial $f$ in the monomial 
order. Let $\langle \LT(I) \rangle$ denote the ideal generated by the 
leading terms of $\LT(f)$ for all $f\in I\setminus \{0\}$. Then it is is well 
known that there is an isomorphism between $\RR[\xb] / I$ and the space 
$\CS_I: = \sspan \{x_1^k x_2^j: x_1^k x_2^j \notin \langle \LT(I) \rangle $
(\cite[Chapt. 5]{CLS}).  The codimension of the ideal is defined by 
$\codim(I): = \dim (\RR[\xb]/I ).$ 

For an ideal $I$ of $\RR[\xb]$, we denote by $V = V(I)$ its real affine 
variety. We consider the case of a zero dimensional variety, that is, when $V$
is a finite set of distinct points in $\RR^2$. In this case, it is well-known 
that $|V| \le \codim I$. The following result is proved in \cite{X00}. 

\begin{prop} \label{prop1}
Let $I$ be a polynomial ideal in $\RR[\xb]$ with finite codimension and 
let $V$ be its affine variety. If $|V| = \codim (I)$ then there is a unique
interpolation polynomial in $\CS_I$ based on the points in the variety. 
\end{prop}

In this case $I$ coincides with the polynomial ideal $I(V)$, which contains
all polynomials in $\RR[\xb]$ that vanish on $V$. An especially interesting
case is when the ideal is generated by a sequence of quasi-orthogonal 
polynomials. Let $d\mu$ be a nonnegative measure with finite moments
on a subset of $\RR^2$. A polynomial $P \in \RR[\xb]$ is said to be an
orthogonal  polynomial with respect to $d\mu$ if 
$$
 \int_{\RR^2} P(\xb) Q(\xb) d\mu(\xb) = 0, \qquad \forall Q\in \RR[\xb],
         \quad \deg Q < \deg P;
$$
it is called a $(2n-1)$-orthogonal polynomial of degree $n$ if the above 
integral is zero for all $Q \in \RR[\xb]$ such that $\deg P + \deg Q \le 
2n-1$.  In particular, if $P$ is of degree $n+1$ and $(2n-1)$ orthogonal, 
then it is orthogonal with respect to all polynomials of degree $n-2$ or 
lower.  We need the following result (cf. \cite{X00}). 

\begin{prop} \label{prop2}
Let $I$ and $V$ be as in the previous proposition. Assume that $I$ is
generated by $(2n-1)$-orthogonal polynomials. If $\codim I = |V|$ then 
there is a cubature formula 
$$
 \int_{\RR^2} f(\xb) d\mu(\xb) = \sum_{\xb \in V}  \lambda_\xb f(\xb) 
$$
of degree $2n-1$; that is, the above formula holds for all $f \in 
\Pi_{2n-1}^2$. 
\end{prop}

The cubature formula in this proposition can be obtained by integrating 
the interpolation polynomial in Proposition \ref{prop1}.  We will apply
this result for $d \mu = W(\xb) d\xb$, where $W$ is the
Chebyshev weight function on $[-1,1]^2$, 
\begin{equation} \label{ChebyW}
   W(\xb) =  \frac{1}{\pi^2} \frac{1}{\sqrt{1-x_1^2}\sqrt{1-x_2^2}}, \qquad   
   -1 < x_1,x_2 < 1.
\end{equation}
For such a weight function, the cubature formula of degree $2n-1$ exists 
only if $|V|$ satisfies M\"oller's lower bound $|V| \ge \dim \Pi_{n-1}^2 
+ [n/2]$. Formulas that attain this lower bound exist (\cite{MP, X94}),
whose knots are the interpolation points studied in \cite{X96}. 

Recall that the Chebyshev polynomials $T_n$ and $U_n$  of the first and 
the second kind are given by $T_n (x) = \cos n \theta$ and $U_n(x) = 
\frac{\sin (n+1)\theta}{\sin \theta}$, $x = \cos \theta$, respectively.  The
polynomials $T_n(x)$ are orthogonal with respect to $1/\sqrt{1-x^2}$. 
Consequently, it is  easy to verify that the polynomials 
\begin{equation} \label{ChebyT}
  T_k(x_1) T_{n-k}(x_2), \qquad 0 \le k \le n, 
\end{equation}
are mutually orthogonal polynomials of degree $n$ with respect to the
weight function $W(x_1,x_2)$ on $[-1,1]^2$.

\subsection{Polynomial ideals and Padua points}

The Padua points \eqref{padua1} were introduced in \cite{CDMV05}  (apart 
from a misprint that $n-1$ should be replaced by $n+1$) based on a 
heuristic argument that we now describe. A set of interpolation points on 
$[-1,1]$ appeared early in \cite{MP}, which are the knots of a cubature formula
of degree $2n-2$ with respect to the weight function 
$\sqrt{1-x_1^2}\sqrt{1-x_2^2}$
on $[-1,1]$. These points are exactly those Padua points \eqref{padua1} that
are inside $(-1,1)^d$ and they are the common zeros of the polynomials 
\begin{equation}\label{MPpoly}
 R_j^{n-1}(x_1, x_2):=
   U_j(x_1) U_{n-j-1}(x_2) +U_{n-j-2}(x_1)U_j(x_2), \qquad 0 \le j \le n-1.
\end{equation}
For $n \ge 1$, there are exactly $\dim \Pi_{n-1}^2$ many such points,
which are not near optimal (\cite{BCDMV05}), however. These points lie 
on various horizontal and vertical lines. Adding proper boundary  points 
of the lines leads to a set of points whose cardinality is exactly 
$\dim \Pi_n^2$, which is the set of Padua points \eqref{padua1}. 

As shown in \cite{BCDMVX} the Padua points can be characterized as
self-intersection points (interior points) and the boundary contact points 
of the algebraic  curve $T_n(x_1)-T_{n+1}(x_2)=0$, $(x_1,x_2)\in [-1,1]^2$ 
(the {\it generating curve}), where $T_n$ denotes the $n$-th Chebyshev 
polynomial of the first kind.  Actually the generating curve used in 
\cite{BCDMVX} is $T_{n+1}(x_1)-T_n(x_2)=0$. Depending on the orientation 
of the $x_1$ and $x_2$ directions, there are four families of Padua points, 
which are the two mentioned above and two others whose generating 
curves are $T_n(x_1)-T_{n+1}(x_2)=0$ and  $T_{n+1}(x_1)+T_n(x_2)=0$, 
respectively. We restrict to the family \eqref{padua1} in this paper. 

\begin{thm}
Let $Q_k^{n+1} \in \Pi_{n+1}^2$ be defined by 
\begin{equation} \label{Q0}
Q_0^{n+1}(x_1,x_2)=T_{n+1}(x_1)-T_{n-1}(x_1), 
\end{equation}
and for $ 1 \le k \le n+1$, 
\begin{equation} \label{Qk}
  Q_k^{n+1}(x_1,x_2) =T_{n-k+1}(x_1)T_k(x_2)+T_{n-k+1}(x_2)T_{k-1}(x_1). 
\end{equation}
Then the set $\Pad_n$ in \eqref{padua1} is the variety of the ideal 
$I = \langle Q_0^{n+1}, Q_1^{n+1}, \ldots, Q_n^{n+1}\rangle$.
Furthermore, $\codim (I) = |\Pad_n| = \dim \Pi_n^2$. 
\end{thm}

\begin{proof}
That the polynomials $Q_k^{n+1}$ vanish on the Padua points \eqref{padua1}
can be easily verified upon using the following representation
$$
Q_0^{n+1}(x_1,x_2)= -2(1-x_1^2)U_{n-1}(x_1) = - 2 \sin \theta \sin n \theta,  
$$
and, for $1\le k \le n+1$, 
\begin{align*}
Q_k^{n+1}(x_1,x_2) = & \cos{k\theta}(\cos{(n+1)\phi}\cos{k\phi}+
   \sin{(n+1)\phi}\sin{k\phi}) \\
 & +\cos{k\phi}(\cos{n\theta}\cos{k\theta} +\sin{n\theta}\sin{k\theta}), 
\end{align*}
where $x_1 = \cos \theta$ and  $x_2 = \cos \phi$. Furthermore, the definition
of $Q_k^{n+1}$ shows easily that the set $\{\LT(Q_0^{n+1}), \ldots, 
\LT(Q_{n+1}^{n+1})\}$ is exactly the set of monomials of degree $n+1$ in
$\RR[\xb]$. Hence, it follows that $\dim \RR[\xb]/I \le  \dim\Pi_n^d$. Recall
that $\codim (I) = \dim \RR[\xb]/I \ge |V|$ and $|V| \ge |\Pad_n| = 
\dim\Pi_n^d$, we conclude that  $\codim (I) = |\Pad_n|$. 
\end{proof}

The basis \eqref{Q0} and \eqref{Qk} of the ideal was identified by using the
fact that the interior points of $\Pad_n$ are the common zeros of 
\eqref{MPpoly}, which shows that the ideal contains the polynomials 
$(1-x_1)^2(1-x_2)^2 R_j^{n-1}$, $0 \le j \le n-1$, as well as two 
polynomials that are actually one variable, 
one being $(1-x_1^2) U_{n-1}(x_1)$, which gives rise to \eqref{Q0}, and the 
other being a polynomial in $x_2$,  which vanishes on $\Pad_n$.  From these 
polynomials we worked out the $Q_k^{n+1}$ whose degree is lower and they in 
fact form the basis of the ideal.  Once the basis is identified, it is easier 
to verify it directly as we did in the proof.

Since the product Chebyshev polynomials \eqref{ChebyT} are orthogonal with
respect to $W$ in \eqref{ChebyW}, it follows readily that $Q_k^{n+1}(x_1,x_2)$
are orthogonal to polynomials in $\Pi_{n-2}^2$ with respect to $W$ 
on $[-1,1]^2$. In particular, they are $(2n-1)$-orthogonal polynomials.
Hence, as a consequence of the theorem and Propositions \ref{prop1} 
and \ref{prop2}, we have the following: 

\begin{cor} \label{cor1}
The set of Padua points $\Pad_n$ is unisolvent for $\Pi_n^d$. Furthermore, 
there is a cubature formula of degree $2n-1$ based on the Padua points in 
$\Pad_n$. 
\end{cor}

We denote the unique interpolation polynomial based on $\Pad_n$ by  
$\CL_n f$, which can be written as 
\begin{equation}\label{Lnf}
 \CL_n f(\xb) = \sum_{\xb_{k,j} \in \Pad_n} f(\xb_{k,j}) \ell_{k,j}(\xb), 
      \qquad \xb = (x_1,x_2), 
\end{equation}
where $\ell_{k,j} (\xb)$ are the fundamental interpolation polynomials 
uniquely determined by 
\begin{equation}\label{ell}
   \ell_{k,j} (\xb_{k',j'})  = \delta_{k,k'}\delta_{j,j'}, \quad
               \hbox{$\forall  \xb_{k,j} \in \Pad_n$, \quad and} \quad   \ell_{k,j} \in \Pi_n^d. 
\end{equation}
In the following section we derive an explicit formula for the fundamental 
polynomials. 

\section{Construction of the Lagrange interpolation polynomials}
\setcounter{equation}{0}

To derive an explicit formula for the Lagrange interpolation formula, we will
use orthogonal polynomials and follow the strategy in \cite{X94,X96}. We note
that the method in \cite{X94} works for the case that $|V| = \dim \Pi_{n-1}^2+
\sigma$ with $\sigma \le n$. In our case $\sigma = n+1$, so that the general
theory there does not apply. 

Let $\CV_n^d$ denote the space of orthogonal polynomials of degree $n$ with 
respect to \eqref{ChebyW} on $[-1,1]^2$. An orthonormal basis for $\CV_n^d$ 
is given by 
$$
P_{n}^{k}(x_1,x_2):=\tilde{T}_{n-k}(x_1) \tilde{T}_k(x_2),  \qquad 0\leq k\leq n. 
$$
where $\wt T_0(x) =  1$ and $\tilde{T}_k(x)=\sqrt{2} T_k(x)$ for $k \ge 1$. 
We introduce the notation 
$$%\begin{equation} \label{basis}
         \PP_n=\left[ P_0^{n},P_1^{n},\dots, P_{n}^{n} \right]
$$%\end{equation}
and treat it both as a set and as a column vector. The reproducing kernel of
the space $\Pi_n^d$ in $L^2(W,[-1,1]^2)$ is defined by
\begin{equation} \label{RK}
   \Kb_n({\xb}, {\yb})=\sum_{k=0}^n [\PP_k (\xb)]^t \PP_k(\yb) 
       = \sum_{k=0}^n \sum_{j=0}^k P_j^k(\xb)P_j^k(\yb). 
\end{equation}
There is a Christoffel-Darboux formula (cf. \cite{DX, X94}) which states that 
\begin{equation} \label{CD}
\Kb_n({\bf x\/},{\bf y\/})= \frac{\left[A_{n,i} 
\PP_{n+1}({\bf x\/})\right]^t\PP_{n}({\bf y\/})-
\left[ A_{n,i} \PP_{n+1}({\bf y\/})\right]^t\PP_{n}({\bf x\/})}{x_i-y_i}, \quad i=1,2\;,
\end{equation}
where ${\bf x\/}=(x_1,x_2)$, ${\bf y\/}=(y_1,y_2)$, and $A_{n,i}$ are matrices
defined by 
$$%\begin{equation} \label{Ai}
A_{n,1}=\frac{1}{2}\,\left(\begin{array} {ccccc}
1 &  & \bigcirc & 0 & 0\\
 & \ddots &  & \vdots & \vdots\\
\bigcirc &  & 1 & 0 & 0\\
0 & \dots & 0 & \sqrt{2} & 0\\
\end{array} \right)\;,\;\;
A_{n,2}=\frac{1}{2}\,\left(\begin{array} {ccccc}
0 & \sqrt{2} & 0 & \dots & 0\\
0 & 0 & 1 &  & \bigcirc\\
\vdots & \vdots &  & \ddots & \\
0 & 0 & \bigcirc & & 1\\
\end{array} \right)\;.
$$%\end{equation}
The fundamental interpolation polynomials are given in terms of 
$\Kb_n(\xb,\yb)$. To this end, we will try to express  $\Kb_n(\xb,\yb)$ 
in terms  of the polynomials in the ideal $I$ of Theorem 2.3. 

Recall that the Padua points are the common zeros of polynomials $Q_k^{n+1}$ 
defined in \eqref{Q0}and \eqref{Qk}. We also denote 
\begin{equation} \label{VectQ}
   \QQ_{n+1}:=\left[\sqrt{2} Q_0^{n+1},2 Q_1^{n+1},\dots,2 Q_n^{n+1},
     \sqrt{2} Q_{n+1}^{n+1} \right]. 
\end{equation}
The definition of $Q_k^{n+1}$ shows that we have the relation 
\begin{equation} \label{3term}
   \QQ_{n+1}=\PP_{n+1}+\Gamma_1\PP_n+\Gamma_2\PP_{n-1},
\end{equation}
where $\Gamma_1$ and $\Gamma_2$ are matrices defined by 
$$%\begin{equation} \label{gammai}
\Gamma_1=\left(\begin{array} {cccc}
0 & \dots & 0 & 0\\
0 & \dots & 0 & \sqrt{2}\\
\bigcirc &  & 1 & 0\\
 & \ddots &  & \vdots\\
1 &  & \bigcirc & 0\\
\end{array} \right)\qquad \hbox{and} \qquad
\Gamma_2=\left(\begin{array} {cc}
-1 & \bigcirc\\ & \\\bigcirc & \bigcirc
\end{array} \right). 
$$%\end{equation}

By using the representation (\ref{3term}) we can rewrite (\ref{CD}) as 
\begin{align} \label{Kn2}
(x_i-y_i)\Kb_n({\bf x\/},{\bf y\/})= & \left[ 
\QQ_{n+1}({\bf x\/})-\Gamma_1\PP_n({\bf x\/})
-\Gamma_2\PP_{n-1}({\bf x\/})\right]^t A_{n,i}^t\PP_n({\bf y\/}) \\
& -\PP_n^t({\bf x\/})A_{n,i}\left[ 
  \QQ_{n+1}({\bf y\/})-\Gamma_1\PP_n({\bf y\/})
  -\Gamma_2\PP_{n-1}({\bf y\/})\right]  \notag \\  
 =  & \Sb_{1,i}({\bf x\/},{\bf y\/})+\Sb_{2,i}({\bf x\/},{\bf y\/})
         +\Sb_{3,i}({\bf x\/},{\bf y\/})\;,\;\;i=1,2 \notag
\end{align}
where
\begin{align*}
\Sb_{1,i}({\bf x\/},{\bf y\/}) & = \QQ^t_{n+1}({\bf x\/})A_{n,i}^t
\PP_n({\bf y\/})-\PP_n^t({\bf x\/})A_{n,i}\QQ_{n+1}({\bf y\/}), \\ 
\Sb_{2,i}({\bf x\/},{\bf y\/}) & = \PP_n^t({\bf x\/})
\left(A_{n,i}\Gamma_1-\Gamma^t_1A_{n,i}^t\right)\PP_n({\bf y\/}), \\
\Sb_{3,i}({\bf x\/},{\bf y\/}) & = 
\PP_n^t({\bf x\/})A_{n,i}\Gamma_2\PP_{n-1}({\bf y\/})
-\PP_{n-1}^t({\bf x\/})\Gamma^t_2A_{n,i}^t\PP_n({\bf y\/})\;.
\end{align*}
If both $\xb$ and $\yb$ are in $\Pad_n$, then $\Sb_{1,i}$ will be zero. We now
work out the other two terms.
First, observe that $A_{n,1}\Gamma_1$ is a symmetrix matrix, 
$$
A_{n,1}\Gamma_1=\left(\begin{array} {ccccc}
0 & & \dots &  & 0\\
\bigcirc &  &  &  & \sqrt{2}\\
  &  &  & 1 & \\
 &  & \ddots &  & \\
  & 1 &  &  & \\ 
  \sqrt{2} &  &  &  & \bigcirc\\
0 & & \dots &  & 0
\end{array} \right), \quad \hbox{and}\quad 
A_{n,1}\Gamma_2=\frac{1}{2}\,
\left(\begin{array} {cc}
- 1 & \bigcirc\\ & \\
\bigcirc & \bigcirc
\end{array} \right)\; .
$$
It follows that $\Sb_{2,1}(\xb,\yb) \equiv 0$ and 
\begin{align}\label{S31}
\Sb_{3,1}({\bf x\/},{\bf y\/})& = 
 -\frac{1}{2}\,P_0^{n}({\bf x\/})P_0^{n-1}({\bf y\/})
    + \frac{1}{2}\,P_0^{n}({\bf y\/})P_0^{n-1}({\bf x\/}) \\
     &  = - T_n(x_1)T_{n-1}(y_1) +T_n(y_1)T_{n-1}(x_1). \notag
\end{align}
We also have $A_{n,2}\Gamma_2=0$, which entails that 
$\Sb_{3,2}(\xb,\yb) \equiv 0$. Finally, we have
$$
A_{n,2}\Gamma_1-\Gamma_1^tA_{n,2}^t=
\frac{1}{2}\,\left(\begin{array} {ccc}
\bigcirc &  & 1\\ & & \\
 &  \bigcirc & \\  & & \\ -1 & & \bigcirc
\end{array} \right)\;,
$$
from which we get 
\begin{equation} \label{S22}
\Sb_{2,2}({\bf x\/},{\bf y\/})=   T_n(x_1)T_{n}(y_2)  - T_n(x_2)T_{n}(y_1). 
\end{equation}
 
The terms $\Sb_{3,1}$ and $\Sb_{2,2}$ do not vanish on the Padua points. We
try to make them part of the left hand side of \eqref{Kn2} by seeking a 
polynomial $h_n({\bf x\/},{\bf y\/})$ such that 
\begin{align} \label{vanish}
\begin{split}
& \Sb_{3,1}({\bf x\/},{\bf y\/}) - (x_1-y_1)h_n({\bf x\/},{\bf y\/}) =0,
 \qquad {\bf x\/},{\bf y\/}\in \Pad_n \;, \\
& \Sb_{2,2}({\bf x\/},{\bf y\/}) - (x_2-y_2)h_n({\bf x\/},{\bf y\/})=0, \qquad 
     {\bf x\/},{\bf y\/}\in \Pad_n. 
\end{split}
\end{align}
It is easy to check that 
$$%\begin{equation} \label{kn}
h_n({\bf x\/},{\bf y\/}):=T_n(x_1)T_n(y_1)  
$$%\end{equation}
satisfies the equations. This can be verified directly using the three term
relation $2 t T_n(t) = T_{n+1}(t) + T_{n-1}(t)$ as follows. By \eqref{S31},
\begin{align*}
&  \Sb_{3,1}({\bf x\/},{\bf y\/})- (x_1-y_1)h_n({\bf x\/},{\bf y\/}) \\
\qquad & = \frac{1}{2}T_n(x_1)(T_{n+1}(y_1) -T_{n-1}(y_1))
   - \frac{1}{2}T_n(y_1)(T_{n+1}(x_1)-T_{n-1}(x_1)) \\
\qquad & = \frac{1}{2}T_n(x_1)Q_0^{n+1}({\bf y\/}) - \frac{1}{2}T_n(y_1)
           Q_0^{n+1}({\bf x\/})\;, 
\end{align*}
using the definition of $Q_0^{n+1}$ at (\ref{Q0}). Similarly, by \eqref{S22},
\begin{align*}
 & \Sb_{2,2}({\bf x\/},{\bf y\/}) - (x_2-y_2)h_n({\bf x\/},{\bf y\/}) \\
  \qquad & = - T_n(y_1)(x_2T_n(x_1)+T_n(x_2))
+ T_n(x_1)(y_2T_n(y_1)+T_n(y_2)) \\
  \qquad & = -  T_n(y_1)Q_1^{n+1}({\bf x\/}) + T_n(x_1)Q_1^{n+1} ({\bf y\/}),
\end{align*}
using the definition of $Q_1^{n+1}$ at \eqref{Qk}.  Since $Q_k^{n+1}(\xb)$ 
vanishes on the Padua points, both equations in \eqref{vanish} are satisfied. 
Consequently, we have proved the following proposition: 

\begin{prop} \label{prop3.1}
For $\xb =(x_1,x_2)$ and $\yb = (y_1,y_2)$, define 
\begin{equation} \label{K*}
    \Kb^*_n(\xb,\yb):= \Kb_n(\xb,\yb) - T_n(x_1)T_n(y_1).
\end{equation}
Then 
\begin{align*}
(x_1-y_1) \Kb^*_n(\xb,\yb)= & \, \QQ^t_{n+1}({\bf x\/})A_{n,1}^t 
\PP_n({\bf y\/})-\PP_n^t({\bf x\/})A_{n,1}\QQ_{n+1}({\bf y\/}) \\
    & +\tfrac{1}{2}T_n(x_1)Q_0^{n+1}({\bf y\/})-\tfrac{1}{2}T_n(y_1)
           Q_0^{n+1}({\bf x\/})\;, \\
(x_2-y_2) \Kb^*_n(\xb,\yb)= & \, \QQ^t_{n+1}({\bf x\/})A_{n,2}^t 
       \PP_n({\bf y\/})-\PP_n^t({\bf x\/})A_{n,2}\QQ_{n+1}({\bf y\/}) \\
    &   - T_n(y_1)Q_1^{n+1}({\bf x\/}) + T_n(x_1)
           Q_1^{n+1}({\bf y\/}). 
\end{align*}
In particular, $\Kb^*(\xb,\yb) =0$ if $\xb, \yb \in \Pad_n$ and $\xb \ne \yb$. 
\end{prop}

Since $P_n^n(x_1,x_2) = T_n(x_1)$, we evidently have $K_n^*(\xb,\xb) 
\ge 1 >0$. Consequently, a compact formula for the fundamental interpolation
polynomials follows immediately from the above proposition.

\begin{thm}
The fundamental interpolation polynomials $\ell_{k,j}$ in \eqref{Lnf} 
associated with the Padua points $\xb_{k,j}= (\xi_k,\eta_j)$ are given by
\begin{equation} \label{compact}
  \ell_{k,j}(\xb) = \frac{\Kb_n^*(\xb,\xb_{k,j})}{\Kb_n^*(\xb_{k,j},\xb_{k,j})}
           =  \frac{\Kb_n(\xb,\xb_{k,j}) - T_n(x_1)T_n(\xi_k)}
                {\Kb_n(\xb_{k,j},\xb_{k,j}) - [T_n(\xi_k)]^2 }. 
\end{equation}
\end{thm}

In fact, using the representation of $\Kb_n^*$ in Proposition \ref{prop3.1}, 
the conditions \eqref{ell} can be verified readily, which proves the theorem. 

We can say more about the formula \eqref{compact}. In fact, a compact formula 
for the reproducing kernel $\Kb_n$ appeared in \cite{X95}, which states that 
\begin{align} \label{compact95}
   \Kb_n(\xb,\yb) = & D_n(\theta_1+\phi_1,\theta_2+ \phi_2) + 
   D_n(\theta_1+\phi_1,\theta_2- \phi_2) \\
     & +D_n(\theta_1-\phi_1,\theta_2+ \phi_2)
        +D_n(\theta_1-\phi_1,\theta_2- \phi_2), \notag
\end{align}
where $\xb = (\cos\theta_1,\cos \theta_2)$, $\yb = (\cos\phi_1,\cos \phi_2)$, 
and
$$
D_n(\alpha,\beta) = \frac{1}{2}\frac{\cos ((n+\frac{1}{2}) \alpha)
\cos \frac{\alpha}{2} -\cos ((n+\frac{1}{2}) \beta)\cos \frac{\beta}{2}} 
{\cos \alpha - \cos \beta}.
$$
Hence, the formula \eqref{compact} provides a compact formula for the 
fundamental interpolation polynomials $\ell_{k,j}$.  Furthermore, the 
explicit formula \eqref{compact95} allows us to derive an explicit formula 
for the denominator $\Kb_n^*(\xb_{k,j},\xb_{k,j})$. The Padua points 
\eqref{padua1} naturally divide into three groups, 
$$
\Pad_n = A_v \cup A_b \cup A_{in},
$$ 
where $A_v$ consists of the two vertex points $((-1)^k, (-1)^{n-1})$ for 
$k=0, n$, $A_b$ consists of the other points on the boundary of $[-1,1]^2$, 
and $A_{in}$ consists of interior points inside $(-1,1)^2$. 

\begin{prop} For $\xb_{k,j} \in \Pad_n$, 
$$
     \Kb_n^*(\xb_{k,j},\xb_{k,j}) = n(n+1) \begin{cases} 
            \frac{1}{2},       & \hbox{if   $\xb_{k,j} \in A_{v}$ } \\ 
            1,  & \hbox{if $\xb_{k,j} \in A_{b}$  } \\   
            2,  & \hbox{if  $ \xb_{k,j} \in A_{in}$}.  \end{cases}
$$
\end{prop} 

The proof can be derived from the explicit formula of the Padua points 
\eqref{padua1} and the compact formula \eqref{compact95} by a tedious 
verification. We refer to \cite{BCDMVX} for a proof using the generating curve.

Since integration of the Lagrange interpolation polynomial gives the quadrature
formula, as a corollary of the above proposition and Corollary \ref{cor1} we
have proved the following:

\begin{prop}
A cubature formula of degree $2n-1$ based on the Padua points is given by
$$
 \frac{1}{\pi^2}  \int_{[-1,1]^2} f(x_1,x_2) 
       \frac{dx_1dx_2}{\sqrt{1-x_1^2}\sqrt{1-x_2^2}}
     = \sum_{\xb_{k,j} \in \Pad_n} w_{k,j} f(\xb_{k,j})
$$
where $w_{k,j}= 1/ \Kb_n^*(\xb_{k,j},\xb_{k,j})$. 
\end{prop}

We note that this cubature formula is not a minimal cubature formula, since 
there are cubature formulas of the same degree with fewer number of nodes.

\section{Convergence of the Lagrange interpolation polynomials}
\setcounter{equation}{0}

In \cite{BCDMVX} it is shown that the Lagrange interpolation $\CL_nf$ in 
\eqref{Lnf} has Lebesgue constant $\CO((\log n)^2)$; that is, as an operator
from $C([-1,1]^2)$ to itself, its operator norm in the uniform norm is 
$\CO( (\log n)^2)$. This settles the problem of uniform convergence of 
$\CL_n f$. In this section we prove that $\CL_n f$ converges in $L^p$ norm. 
For the Chebyshev weight function $W$ defined in \eqref{ChebyW}, we 
define $L^p(W)$ as the space of Lebesgue measurable functions for
which the norm 
$$
\|f\|_{W,p} : = \left (\int_{[-1,1]^2} |f(x_1,x_2)|^p W(x_1,x_2)
     dx_1dx_2\right)^{1/p}
$$
is finite. We keep this notation also for $0 < p<1$, even though it is no
longer a norm for $p$ in that range. The main result in this section is

\begin{thm}\label{thm4.1} 
 Let $\CL_nf$ be the Lagrange interpolation polynomial \eqref{Lnf} based on 
 the Padua points. Let $0 < p < \infty$. Then
 $$
\lim_{n\to \infty} \|L_n f - f\|_{W,p} = 0, \qquad \forall  f\in C([-1,1]^2). 
 $$
In fact, let $E_n(f)_\infty $ be the error of the best approximation of $f$ by 
polynomials from $\Pi_n^2$ in the uniform norm; then 
$$
\|L_n f - f\|_{W,p}  \le c_p E_n (f)_\infty  \qquad \forall  f\in C([-1,1]^2). 
$$
\end{thm} 
 
The proof follows the approach in \cite{X96}, where the mean convergence 
of another family of interpolation polynomials is proved. We shall be brief 
whenever the same proof carries over. First we need a lemma on the Fourier 
partial sum, $S_n f$, defined by 
$$
   S_n f (\xb) = \sum_{k=0}^n \sum_{j=0}^k a_j^k(f) P_j^k(\xb)  
         =  \frac{1}{\pi^2} \int_{[-1,1]^2} \Kb_nf(\xb,\yb) W(\yb) d\yb
$$
where $a_j^k(f)$ is the Fourier coefficient of $f$ with respect to the orthonormal
basis $\{P_j^k\}$ in $L^2(W)$. In the following we let $c_p$ denote a 
generic constant that depends on $p$ only, its value may be different from 
line to line.  

\begin{lem}  \label{lem4.1}
Let $1<p<\infty$. Then 
\begin{equation} \label{4.1}
      \|S_n f\|_{W,p} \le c_p \|f\|_{W,p}, \qquad \forall f \in L^p(W).
\end{equation}
\end{lem}

This is \cite[Lemma 3.4]{X96}. We will need another lemma whose proof follows 
almost verbatim from that of \cite[Lemma 3.5]{X96}.  In the following we write
$$
     N = | \Pad_n| = (n+1)(n+2)/2.
$$

\begin{lem}  \label{lem4.2}
Let $1 \le p \le \infty$. Let $\xb_{k,j}$ be the Padua points. Then  
\begin{equation} \label{4.2}
  \frac{1}{N} \sum_{\xb_{k,j} \in \Pad_n} |P(\xb_{k,j})|^p \le 
   c_p \int_{[-1,1]^p} |P(\xb)|^p W(\xb) d\xb, \qquad \forall P\in \Pi_n^2. 
\end{equation}
\end{lem}  
 
The main tool in the proof of Theorem 4.1 is a converse inequality of 
\eqref{4.2}, which we state below and give a complete proof,  even though 
the proof is similar to that of \cite[Theorem 3.3]{X96}.

\begin{prop}
Let $1 < p< \infty$. Let $\xb_{k,j}$ be the Padua points. Then 
$$
    \int_{[-1,1]^2} |P(\xb)|^p W(\xb) d\xb \le c_p \frac{1}{N}
     \sum_{\xb_{k,j} \in \Pad_n} |P(\xb_{k,j}|^p, \qquad \forall P \in \Pi_n^d.
$$
\end{prop}

\begin{proof}
Let $P \in \Pi_n^d$. For $p > 1$, we have 
\begin{align} \label{norm}
  \|P\|_{W,p} = & \sup_{\|g\|_{W,q} =1} \int_{[-1,1]^2} P(\xb) g(\xb) W(\xb) d\xb \\
                    = & \sup_{\|g\|_{W,q} =1} \int_{[-1,1]^2} P(\xb) S_n g(\xb) W(\xb) d\xb,                 
     \qquad \frac{1}{p} +\frac{1}{q} =1, \notag
\end{align}
where the second equality follows from the orthogonality. We write 
$$
   S_n g = S_{n-1} g + \ab_n^t(g) \PP_n, \quad\hbox{where}\quad
       \ab_n=\int_{[-1,1]^2} g(\xb)\PP_n(\xb) W(\xb)d\xb. 
$$
Since the cubature formula is of degree $2n-1$ and $P S_{n-1}g$ is of degree
$2n-1$, we have 
\begin{align} \label{4.4}
 & \left| \int_{[-1,1]^2} P(\xb) S_{n-1}(\xb) W(\xb) d\xb \right | = 
       \left |\sum_{\xb_{k,j} \in\Pad_n} w_{k,j} P(\xb_{k,j}) 
   S_{n-1}(\xb_{k,j}) \right|    \notag \\
 &  \qquad    \le \left (\sum_{\xb_{k,j} \in \Pad_n} w_{k,j} 
   |P(\xb_{k,j})|^p \right)^{1/p} \left (\sum_{\xb_{k,j} \in\Pad_n} 
     w_{k,j}  |S_{n-1}(\xb_{k,j})|^q \right)^{1/q}. 
\end{align}
By Lemma \ref{lem4.1} and Lemma \ref{lem4.2},  
$$
 \sum_{\xb_{k,j} \in\Pad_n} w_{k,j}  |S_{n-1}(\xb_{k,j})|^q  \le 
c \|S_{n-1} g\|_{W,p} \le c \|g\|_{W,p}= c.
$$
Hence,  using the fact that $w_{k,j} \sim N^{-1}$, it follows readily that 
\begin{equation} \label{Sn-1}
\left| \int_{[-1,1]^2} P(\xb) S_{n-1}(\xb) W(\xb) d\xb \right | \le  c
   \left (\sum_{\xb_{k,j} \in \Pad_n} w_{k,j} |P(\xb_{k,j})|^p \right)^{1/p}. 
\end{equation}
We need to establish the similar inequality for $\ab_n^t \PP_n$ term. Since 
$\CL_n P = P$ as  $P$ is of degree $n$, it follows from \eqref{compact} and
the orthogonality that
\begin{align*}
& \int_{[-1,1]^2} P(\xb) \ab_n^t(g)\PP_n(\xb) W(\xb) d\xb =
  \int_{[-1,1]^2} \CL_nP(\xb) \ab_n^t(g)\PP_n(\xb) W(\xb) d\xb \\
& =  \sum_{\xb_{k,j} \in\Pad_n} w_{k,j} P(\xb_{k,j}) 
           \int_{[-1,1]^2} \Kb_n^* (\xb_{k,j}, \xb) \ab_n^t(g)\PP_n(\xb) W(\xb) d\xb \\
& = \sum_{\xb_{k,j} \in\Pad_n} w_{k,j} P(\xb_{k,j}) 
           \ab_n^t(g)\PP_n(\xb_{k,j}) -  \frac{1}{2} a_n^n(g) \sum_{\xb_{k,j} \in\Pad_n} w_{k,j} 
             P(\xb_{k,j}) \wt T_n(\xi_{k,j}).
\end{align*}
For the first term we can write $\ab^t(g) \PP_n= S_ng - S_{n-1}g$ and apply 
H\"older's inequality as in \eqref{4.4}, so that the proof of \eqref{Sn-1} can be
carried out again. For the second term we use the fact that $|\wt T_n(x)| \le 
\sqrt{2}$ to conclude that 
$$
   |a_n^n(g)| \le \int_{[-1,1]^2} |\wt T_n(\xb) g(\xb)| W(\xb) d\xb 
      \le  \sqrt{2} \|g\|_{W,q} \le \sqrt{2}.   
$$
so that the sum is bounded by
\begin{align*}
  \left| a_n^n(g) \sum_{\xb_{k,j} \in\Pad_n} w_{k,j} 
             P(\xb_{k,j}) \wt T_n(\xi_{k,j})\right| \le \sqrt{2} \left(
             \sum_{\xb_{k,j} \in\Pad_n} w_{k,j} |P(\xb_{k,j})|^p \right)^{1/p}.
\end{align*}
This way we have established the inequality 
$$
  \left| \int_{[-1,1]^2} P(\xb) \ab_n^t(g)\PP_n(\xb) W(\xb) d\xb \right|
    \le c \left (\sum_{\xb_{k,j} \in \Pad_n} w_{k,j} |P(\xb_{k,j})|^p \right)^{1/p}. 
$$
Together with \eqref{Sn-1}, this completes the proof of the proposition.
\end{proof}

As shown in the proof of Theorem 3.1 of \cite{X96}, the proof of Theorem 4.1
follows as an easy consequence of the Proposition 4.4.


\begin{thebibliography}{99}
    
\bibitem{BCDMV05} 
    L.~Bos, M.~Caliari, S.~De Marchi and M.~Vianello,
     \textit{A numerical study of the Xu polynomial interpolation formula in 
     two variables}. Computing {\bf 76\/} (2005), 311--324.

\bibitem{BDM} 
     L.~Bos, S.~De Marchi and M.~Vianello, 
     \textit{On the Lebesgue constant for the Xu interpolation
     formula}. J. Approx. Theory, to appear.

\bibitem{BCDMVX} 
     L.~Bos, M.~Caliari, S.~De Marchi, M.~Vianello and Y.~Xu, 
    \textit{Bivariate Lagrange interpolation at the Padua points: 
    the generating curve approach}. J. Approx. Theory, to appear.

\bibitem{CLS}
    D. Cox, J. Little and D. O'Shea,  
    \textit{Ideals, Varieties, and Algorithms}, 2nd ed., Springer, Berlin, 1997.

\bibitem{DX} 
    C.F.~Dunkl and Y. Xu, 
     \textit{Orthogonal Polynomials of Several Variables},
     Encyclopedia of Mathematics and its Applications, 81.  
     Cambridge University Press, Cambridge, 2001.

\bibitem{CDMV05} 
    M.~Caliari, S.~De Marchi and M.~Vianello,
    \textit{Bivariate polynomial interpolation on the square at new
    nodal sets}. Appl. Math. Comput. {\bf 165\/} (2005), 261--274.

\bibitem{MP}
    C. R. Morrow and T. N. L. Patterson, 
    \textit{Construction of algebraic cubature rules using polynomial
    ideal theory}, SIAM J. Numer. Anal.  {\bf 15} (1978), 953--976. 

\bibitem{X94} 
    Y.~Xu,  \textit{Common zeros of polynomials in several variables 
    and higher dimensional quadrature}, Pitman Research Notes in
     Mathematics Series, Longman, Essex, 1994.
     
\bibitem{X95} 
    Y.~Xu, \textit{Christoffel functions and Fourier series for 
    multivariate orthogonal polynomials}. J. Approx. Theory {\bf 82\/} (1995), 
     205--239.

\bibitem{X96} 
    Y.~Xu, \textit{Lagrange interpolation on Chebyshev points of 
    two variables}. J. Approx. Theory {\bf 87\/} (1996), 220--238.

%\bibitem{X099} Y.~Xu, \textit{Cubature formulae and polynomial 
%ideals}. Adv. in Appl. Math. 23 (1999), 211--233.

\bibitem{X00} Y.~Xu, 
    \textit{Polynomial interpolation in several variables, cubature 
     formulae, and ideals}, Adv. Comput. Math. {\bf 12\/} (2000), 363--376.

\end{thebibliography}
\end{document}